# Companion forms and weight one forms

By Kevin Buzzard and Richard Taylor*

## Introduction

In this paper we prove the following theorem. Let $L/\mathbb{Q}_p$ be a finite extension with ring of integers $\mathcal{O}_L$ and maximal ideal $\lambda$.

THEOREM 1. *Suppose that $p \geq 5$. Suppose also that $\rho : G_\mathbb{Q} \to \mathrm{GL}_2(\mathcal{O}_L)$ is a continuous representation satisfying the following conditions.*

1. *$\rho$ ramifies at only finitely many primes.*

2. *$(\rho \bmod \lambda)$ is modular and absolutely irreducible.*

3. *$\rho$ is unramified at $p$ and $\rho(\mathrm{Frob}_p)$ has eigenvalues $\alpha$ and $\beta$ with distinct reductions modulo $\lambda$.*

*Then there exists a classical weight one eigenform $f = \sum_{n=1}^\infty a_m(f) q^m$ and an embedding of $\mathbb{Q}(a_m(f))$ into $L$ such that for almost all primes $q$, $a_q(f) = \mathrm{tr}\, \rho(\mathrm{Frob}_q)$. In particular $\rho$ has finite image and for any embedding $i : L \hookrightarrow \mathbb{C}$ the Artin L-function $L(i \circ \rho, s)$ is entire.*

We have three motivations for looking at this theorem. Firstly it can be seen as partial confirmation of a conjecture of Fontaine and Mazur which asserts that if $\rho : G_\mathbb{Q} \to \mathrm{GL}_n(\overline{\mathbb{Q}}_p)$ is a continuous irreducible representation ramified at only finitely many primes and such that the image of the inertia group at $p$ is finite, then the image of $\rho$ is finite (see [FM]). Our theorem verifies this conjecture in the case that $n = 2$, $p \geq 5$ and the reduction $(\rho \bmod \lambda)$ is modular, irreducible and takes $\mathrm{Frob}_p$ to an element with distinct eigenvalues. In our opinion the most serious assumption here is that $(\rho \bmod \lambda)$ should be modular, but we remind the reader that if $(\rho \bmod \lambda)$ is odd then a conjecture of Serre (see [S2]) predicts that it is necessarily modular.

*The first author was supported by a Miller Fellowship. He would also like to thank the Harvard University Clay Fund for funding a short visit to Harvard, during which time some of this work was completed. The second author was partially supported by NSF Grant DMS-9702885.

Our second motivation is that our theorem gives a partial answer to the question: When is the weight one specialisation of a Hida family of ordinary eigenforms a classical weight one form? This question was raised by Mazur and Wiles (see [MW]), who showed that the answer is not "always". We expect the answer to be: Whenever the image of inertia at $p$ in the corresponding specialisation of the associated Galois representation is finite. Our theorem gives partial confirmation of this expectation (see Corollary 3).

Our third motivation was that this work forms part of a program, outlined by one of us (RT) in 1992 (first to Wiles and later to other people), to prove the Artin conjecture for certain odd two dimensional icosahedral representations. If our main theorem could be extended to include the case $p = 2$, subject to the additional hypothesis that $\rho$ is odd, then combining it with the main theorem of [ST] we would obtain the Artin conjecture in this case, subject to some local conditions (certainly in the case the representation is unramified at 2, 3 and 5 and the inertial degree of 2 in the $A_5$ extension is $> 2$). The basic arguments of this paper seem to work for $p = 2$, but unfortunately some of our references exclude this prime. We believe the only place this poses a serious problem is the reference to the work of Wiles [W], Taylor-Wiles [TW] and Diamond [D]. Extensions of these results in a suitable way to the case $p = 2$ are currently being investigated by Mark Dickinson for his Harvard Ph.D.

Since writing the original draft of this paper we have discovered that our methods also give partial, but significant, answers both to a question of Gouvea and a question of Serre. The former concerns for which eigenvalues of $U_p$ a $p$-adic modular form can be overconvergent (see [Go, p. 53]). The latter concerns when a mod $p$ weight one eigenform of level prime to $p$ can be lifted into characteristic zero. These applications are described at the end of Sections 2 and 3 respectively.

Outlining the proof of our main theorem, we first use results of Gross [Gr] to find two weight 2 mod $p$ eigenforms $\overline{f}_\alpha$ and $\overline{f}_\beta$ of level $Np$ (where $p \nmid N$) whose associated Galois representations are $(\rho \bmod \lambda)$. They differ however in that the eigenvalues of $U_p$ on $\overline{f}_\alpha$ and $\overline{f}_\beta$ are $\alpha$ and $\beta$. Then we use the results of Wiles [W] as completed by Taylor and Wiles [TW] and extended by Diamond [D] to show that there are $\Lambda$-adic eigenforms (in the sense of Hida [H2]) $F_\alpha$ and $F_\beta$ lifting $\overline{f}_\alpha$ and $\overline{f}_\beta$ respectively, whose associated Galois representations specialise in weight one to $\rho$. Specialising $F_\alpha$ and $F_\beta$ in weight one gives two overconvergent $p$-adic weight one modular forms $f_\alpha$ and $f_\beta$. We let $f = (\alpha f_\alpha - \beta f_\beta)/(\alpha - \beta)$ and $f' = (f_\alpha - f_\beta)/(\alpha - \beta)$. There are two natural projections $\pi_1, \pi_2 : X_0(p) \to X$. We show that, restricted to a certain rigid subspace, $\pi_2^* f' = \pi_1^* h$ for some $h$ defined on a rigid subspace of $X = X_1(N)$ including the central parts of the supersingular discs. Moreover we show that $h$ and $f$ glue together to give a weight one form defined on the whole of $X$. This is the weight one form we are looking for.

It is a pleasure to acknowledge the debt this work owes to the work of Gross [Gr] and to the various works of Coleman on $p$-adic modular forms. We would also like to thank Brian Conrad, Naomi Jochnowitz, Barry Mazur and particularly Robert Coleman for useful conversations.

## 1. $\Lambda$-adic companion forms

Throughout this paper we will fix a prime $p \geq 5$ and an integer $N \geq 5$ which is not divisible by $p$.

We will let $S_k(\Gamma_1(M))$ denote the space of weight $k$ cusp forms on $\Gamma_1(M)$ with rational integer $q$-expansion at infinity. This space comes with an action of $(\mathbb{Z}/M\mathbb{Z})^\times$ (via the diamond operators $x \mapsto \langle x \rangle$), the Hecke operators $T_q$ for any prime $q \nmid M$, and the Hecke operators $U_q$ for any prime $q|M$. For any prime $q \nmid M$ we define $S_q = q^{k-1}\langle q \rangle$ (note that we have not followed our normal convention for the normalisation of $S_q$), and for any $n \in \mathbb{Z}_{\geq 1}$ we define $T(n)$ by the following formulae:

- $T(n_1 n_2) = T(n_1)T(n_2)$ if $n_1$ and $n_2$ are coprime;
- $\sum_{r=0}^{\infty} T(q^r)X^r = (1 - T_q X + S_q X^2)^{-1}$ if $q \nmid M$;
- $T(q^r) = U_q^r$ if $q|M$.

Let $h_k(M)$ denote the $\mathbb{Z}$-algebra generated by all these Hecke operators acting on $S_k(\Gamma_1(M))$. There is a perfect duality

$$\begin{array}{rcl} S_k(\Gamma_1(M)) \times h_k(M) & \longrightarrow & \mathbb{Z} \\ (f, T) & \longmapsto & a_1(f|T), \end{array}$$

where $a_m(g)$ denotes the coefficient of $q^m$ in the $q$-expansion at infinity of $g$.

Whenever it makes sense we will let $e = \lim_{r \to \infty} U_p^{r!}$, the Hida idempotent. Following Hida we set

$$h^0(N) = \varprojlim e(h_2(Np^r) \otimes_\mathbb{Z} \mathbb{Z}_p).$$

Note that the operators $T(n)$ for $n \geq 1$ and $S_q$ for $q$ not dividing $Np$ are compatible with the projection maps and commute with $e$, and so give rise to operators $T(n)$ and $S_q$ in $h^0(N)$. Moreover we have a natural map

$$\langle \ \rangle : \mathbb{Z}_p^\times \times (\mathbb{Z}/N\mathbb{Z})^\times = \varprojlim (\mathbb{Z}/Np^r\mathbb{Z})^\times \to h^0(N)^\times.$$

We let $\Lambda = \mathbb{Z}_p[[(1+p\mathbb{Z}_p)^\times]]$ and $u = (1+p) \in (1+p\mathbb{Z}_p)^\times \subset \Lambda^\times$. Thus $\Lambda \cong \mathbb{Z}_p[[(u-1)]]$. If $\zeta$ is a $p$-power root of unity we let $\psi_\zeta : (1+p\mathbb{Z}_p)^\times \to \mathbb{Z}_p[\zeta]^\times$ denote the continuous homomorphism taking $u$ to $\zeta$. We also let it denote the corresponding homomorphism $\Lambda \to \mathbb{Z}_p[\zeta]$. Note that $h^0(N)$ is a $\Lambda$-module via $(1 + p\mathbb{Z}_p)^\times \to h^0(N)^\times$ by $u \mapsto \langle u \rangle$. Hida (see [H2]) proves that $h^0(N)$ is a finite free $\Lambda$-module and that for any $k \geq 2$ and $\zeta$ a primitive $p^r$ root of unity

$$(h^0(N) \otimes_{\mathbb{Z}_p} \mathbb{Z}_p[\zeta])/(u - (1+p)^{k-2}\zeta)(h^0(N) \otimes_{\mathbb{Z}_p} \mathbb{Z}_p[\zeta])$$

is isomorphic to the maximal quotient of $e(h_k(Np^{r+1}) \otimes_{\mathbb{Z}} \mathbb{Z}_p[\zeta])$ where $\langle x \rangle = \psi_\zeta(x)$ for all $x \in (\mathbb{Z}/Np^{r+1}\mathbb{Z})^\times$ with $x \equiv 1 \bmod Np$. This isomorphism further takes the Hecke operators $T(n)$, $T_q$, $S_q$ and $U_q$ to themselves. For $k = 1$ it is known that in general such an isomorphism does not exist; see for instance [MW].

We set $S^0(N) = \mathrm{Hom}\,_\Lambda(h^0(N), \Lambda)$. Elements $F \in S^0(N)$ are called $\Lambda$-adic forms. There is an injection $S^0(N) \hookrightarrow \Lambda[[q]]$ which takes $F$ to $\sum_{n=1}^\infty F(T(n))q^n$ and which is called the $q$-expansion map. If $\mathcal{L}$ is a finite field extension of the field of fractions of $\Lambda$ then we call $F \in S^0(N) \otimes_\Lambda \mathcal{L} = \mathrm{Hom}\,_\Lambda(h^0(N), \mathcal{L})$ a $\Lambda$-adic eigenform if it is a $\Lambda$-algebra homomorphism $h^0(N) \to \mathcal{L}$. We call two $\Lambda$-adic eigenforms $F_1 \in S^0(N) \otimes_\Lambda \mathcal{L}_1$ and $F_2 \in S^0(N) \otimes_\Lambda \mathcal{L}_2$ equivalent if there are a finite extension $\mathcal{L}_3$ of the field of fractions of $\Lambda$ and embeddings of $\Lambda$-algebras $\mathcal{L}_1 \hookrightarrow \mathcal{L}_3$ and $\mathcal{L}_2 \hookrightarrow \mathcal{L}_3$ which send $F_1$ and $F_2$ to the same element of $S^0(N) \otimes_\Lambda \mathcal{L}_3$. Equivalence classes of $\Lambda$-adic eigenforms are in bijection with height zero primes of $h^0(N)$ via the map which sends an eigenform to its kernel. If $k \in \mathbb{Z}_{\geq 2}$ and $\zeta$ is a primitive $p^r$th root of unity then the map

$$F \longmapsto \sum_{n=1}^\infty (F(T(n)) \bmod (u - (1+p)^{k-2}\zeta))q^n$$

gives rise to an isomorphism of

$$(S^0(N) \otimes_{\mathbb{Z}_p} \mathbb{Z}_p[\zeta])/(u - (1+p)^{k-2}\zeta)(S^0(N) \otimes_{\mathbb{Z}_p} \mathbb{Z}_p[\zeta])$$

with the maximal subspace of $e(S_k(\Gamma_1(Np^{r+1})) \otimes_{\mathbb{Z}} \mathbb{Z}_p[\zeta])$ where $\langle x \rangle = \psi_\zeta(x)$ for all $x \in (\mathbb{Z}/Np^{r+1}\mathbb{Z})^\times$ with $x \equiv 1 \bmod Np$.

We will let $\chi^{\mathrm{cyclo}}$ denote the usual character

$$G_\mathbb{Q} \twoheadrightarrow \mathrm{Gal}\,(\mathbb{Q}(\zeta_{Np^\infty})/\mathbb{Q}) \cong \mathbb{Z}_p^\times \times (\mathbb{Z}/N\mathbb{Z})^\times \hookrightarrow \mathbb{Z}_p[[\mathbb{Z}_p^\times \times (\mathbb{Z}/N\mathbb{Z})^\times]]^\times$$
$$= \Lambda[(\mathbb{Z}/Np\mathbb{Z})^\times]^\times.$$

Also note that the map $q \mapsto S_q$ extends (when we embed the primes not dividing $Np$ diagonally in $\mathbb{Z}_p^\times \times (\mathbb{Z}/N\mathbb{Z})^\times$) to a homomorphism $S : \mathbb{Z}_p^\times \times (\mathbb{Z}/N\mathbb{Z})^\times \to h^0(N)^\times$ which sends $(x_p, x^p)$ to $x_p \langle (x_p, x^p) \rangle$ for $x_p \in \mathbb{Z}_p^\times$ and $x^p \in (\mathbb{Z}/N\mathbb{Z})^\times$. Thus $S = S_p \times S^p$ where $S_p : \mathbb{Z}_p^\times \to h^0(N)^\times$ and $S^p : (\mathbb{Z}/N\mathbb{Z})^\times \to h^0(N)^\times$. Also $S$ extends to a homomorphism $\Lambda[(\mathbb{Z}/Np\mathbb{Z})^\times] \to h^0(N)$.

If $\mathfrak{m}$ is a maximal ideal of $h^0(N)$ we will let $k(\mathfrak{m})$ denote its residue field. There is a unique (up to conjugation) continuous semisimple representation $\overline{\rho}_\mathfrak{m} : G_\mathbb{Q} \to \mathrm{GL}_2(k(\mathfrak{m}))$ such that for all primes $q \nmid Np$ the representation is unramified and $\mathrm{tr}\,\overline{\rho}_\mathfrak{m}(\mathrm{Frob}_q) = T_q$. We call $\mathfrak{m}$ Eisenstein (respectively non-Eisenstein) if $\overline{\rho}_\mathfrak{m}$ is absolutely reducible (resp. absolutely irreducible). If $\mathfrak{m}$ is a non-Eisenstein maximal ideal of $h^0(N)$ then there is a unique (up to conjugation) continuous representation

$$\rho_\mathfrak{m} : G_\mathbb{Q} \longrightarrow \mathrm{GL}_2(h^0(N)_\mathfrak{m})$$

such that for primes $q \nmid Np$, $\rho_\mathfrak{m}$ is unramified at $q$ and $\operatorname{tr} \rho_\mathfrak{m}(\operatorname{Frob}_q) = T_q$. (See [H2].) Moreover $\det \rho_\mathfrak{m} = S \circ \chi^{\mathrm{cyclo}}$ and

$$\rho_\mathfrak{m}|_{G_p} \sim \begin{pmatrix} \phi^{-1}(S \circ \chi^{\mathrm{cyclo}}) & * \\ 0 & \phi \end{pmatrix},$$

where $\phi$ is the unramified character with $\phi(\operatorname{Frob}_p) = U_p$. We call a $\Lambda$-adic eigenform $F \in S^0(N) \otimes_\Lambda \mathcal{L}$ Eisenstein (resp. non-Eisenstein) if the unique maximal ideal of $h^0(N)$ above $\ker F$ is Eisenstein (resp. non-Eisenstein). If $F$ is non-Eisenstein we obtain a continuous representation

$$\rho_F : G_\mathbb{Q} \longrightarrow \mathrm{GL}_2(\mathcal{O}_\mathcal{L})$$

which for all primes $q \nmid Np$ is unramified and satisfies $\operatorname{tr} \rho_F(\operatorname{Frob}_q) = F(T_q)$. Here $\mathcal{O}_\mathcal{L}$ denotes the integral closure of $\Lambda$ in $\mathcal{L}$.

We call two $\Lambda$-adic eigenforms $F$ and $G \in S^0(N) \otimes \mathcal{L}$ *companion forms*, with respect to height one primes $\wp$ and $\wp'$ of $\mathcal{O}_\mathcal{L}$, which do not divide $p$, if we can find embeddings $\mathcal{O}_\mathcal{L}/\wp \hookrightarrow \mathbb{Q}_p^{ac}$ and $\mathcal{O}_\mathcal{L}/\wp' \hookrightarrow \mathbb{Q}_p^{ac}$ such that

1. for all $m \in \mathbb{Z}_{\geq 1}$ not divisible by $p$
$$G(T(m)) \bmod \wp' = F(T(m)S_p(m)^{-1}) \bmod \wp,$$

2. $G(U_p) \bmod \wp' = F(U_p^{-1} S^p(p)) \bmod \wp$.

Note that this is equivalent to

1. $(\rho_G \bmod \wp') \sim (\rho_F \otimes (F \circ S_p \circ \chi^{\mathrm{cyclo}})^{-1} \bmod \wp)$;

2. if $q|N$ then $G(U_q) \bmod \wp' = F(U_q S_p(q)^{-1}) \bmod \wp$.

3. In the case $F \circ S_p \equiv 1 \bmod \wp$ we also require
$$G(U_p) \bmod \wp' = F(U_p^{-1} S^p(p)) \bmod \wp.$$

Note that we will only use this definition in the case where $\overline{\rho}_\mathfrak{m}$ is irreducible and the $(\overline{\rho}_\mathfrak{m}|_{G_p})^{ss}$ do not act as scalars, where $\mathfrak{m}$ is the maximal ideal of $h^0(N)$ containing $\ker F$. It may be that one should modify this definition in other cases. From now on we suppose that $\mathcal{L}$ is fixed as a Galois extension of the fraction field of $\Lambda$ large enough that all $\Lambda$-adic eigenforms of level $N$ are equivalent to ones with values in $\mathcal{L}$.

THEOREM 2. *Suppose $F$ is a $\Lambda$-adic eigenform in $S^0(N) \otimes_\Lambda \mathcal{L}$ ($\mathcal{L}$ as above) and let $\mathfrak{m}$ be the maximal ideal of $h^0(N)$ containing $\ker F$. Suppose that* for all $q|N$, $F(U_q) = 0$. Suppose also that $\overline{\rho}_{\mathfrak{m}}|_{G_{\mathbb{Q}(\sqrt{(-1)^{(p-1)/2}p})}}$ is absolutely irreducible and $(\overline{\rho}_{\mathfrak{m}}|_{G_p})^{ss}$ does not consist of scalar matrices. Let $\wp$ be a height one prime of $\mathcal{O}_{\mathcal{L}}$ not dividing $p$. Then the following are equivalent.

1. $(\rho_F \bmod \wp)|_{G_p} \sim \begin{pmatrix} * & 0 \\ 0 & * \end{pmatrix}$.

2. $F$ has a companion form $G$ with respect to $\wp$ and some second height one prime $\wp'$ of $\mathcal{O}_{\mathcal{L}}$ which does not divide $p$.

Before proving this theorem we make two remarks. First, the theorem is the $\Lambda$-adic analogue of a conjecture of Serre which was proved by Gross in his beautiful paper [Gr]. Second, the condition that $F(U_q) = 0$ for $q|N$ is mostly for simplicity. All we really need assume is that $F(U_q) = 0$ if $q|N$ and $\overline{\rho}_{\mathfrak{m}}$ is unramified at $q$. In any case, if $F$ is an eigenform of level $N$ and if $N' = N \prod_{q|N} q$ then there is an eigenform $F'$ of level $N'$ with $\rho_{F'} = \rho_F$, $F'(U_q) = 0$ for all $q|N$ and $F'(T(n)) = F(T(n))$ for all $n$ coprime to $N$.

We now turn to the proof of the theorem. We first show that condition two implies condition one. We know that

$$(\rho_G \bmod \wp') \sim (\rho_F \otimes (F \circ S_p \circ \chi^{\text{cyclo}})^{-1} \bmod \wp).$$

But we also know that

$$(\rho_G \bmod \wp')|_{G_p} \sim \begin{pmatrix} \phi(F \circ S_p \circ \chi^{\text{cyclo}})^{-1} & * \\ 0 & \phi^{-1}(F \circ S^p \circ \chi^{\text{cyclo}}) \end{pmatrix} \bmod \wp$$

where $\phi$ is unramified with $\phi(\text{Frob}_p) = F(U_p)$. On the other hand

$$(\rho_F \otimes (F \circ S_p \circ \chi^{\text{cyclo}})^{-1} \bmod \wp)|_{G_p}$$
$$\sim \begin{pmatrix} \phi^{-1}(F \circ S^p \circ \chi^{\text{cyclo}}) & * \\ 0 & \phi(F \circ S_p \circ \chi^{\text{cyclo}})^{-1} \end{pmatrix} \bmod \wp.$$

Thus we must have

$$(\rho_G \bmod \wp)|_{G_p} \sim \begin{pmatrix} * & 0 \\ 0 & * \end{pmatrix},$$

as desired.

The reverse implication is much deeper. First note that by [Gr] (see also [CV] where the unproved hypotheses of [Gr] are removed) there is a maximal ideal $\mathfrak{n}$ of $h^0(N)$ such that

$$\overline{\rho}_{\mathfrak{n}} \sim \overline{\rho}_{\mathfrak{m}} \otimes (S_p \circ \chi^{\text{cyclo}} \bmod \mathfrak{m})^{-1};$$

$(U_p \bmod \mathfrak{n}) = (U_p^{-1} S^p(p) \bmod \mathfrak{m})$; and $U_q \equiv 0 \bmod \mathfrak{n}$ for $q|N$ (after some embedding of $k(\mathfrak{n})$ into $k(\mathfrak{m})$). Consider deformations $\rho$ of $\overline{\rho}_{\mathfrak{n}}$ to complete noetherian local $W(k(\mathfrak{m}))$-algebras with residue field $k(\mathfrak{m})$ satisfying

- $\rho$ is unramified outside $Np$,

- $\rho|_{G_p} \sim \begin{pmatrix} * & * \\ 0 & \phi \end{pmatrix}$ with $\phi$ unramified.

Let $\rho^{\mathrm{univ}} : G_\mathbb{Q} \to \mathrm{GL}_2(R)$ be the universal such deformation. Then by Theorem 1.1 of [D] there is a surjection $h^0(N) \twoheadrightarrow R$ taking

- $T_q$ to $\mathrm{tr}\,\rho^{\mathrm{univ}}(\mathrm{Frob}_q)$ for $q \nmid Np$,

- $S_q$ to $\det \rho^{\mathrm{univ}}(\mathrm{Frob}_q)$ for $q \nmid Np$,

- $U_q$ to 0 for $q|N$,

- $U_p$ to $\phi^{\mathrm{univ}}(\mathrm{Frob}_p)$.

(Note that if $q|N$ then $q^{a+b}|N$ where $a$ denotes the conductor of $\overline{\rho}_\mathfrak{n}$ at $q$ and $b = \dim \overline{\rho}_\mathfrak{n}^{I_q}$. This is because $U_q \in \mathfrak{n}$.) However $(\rho_F \otimes (F \circ S_p \circ \chi^{\mathrm{cyclo}})^{-1} \mod \wp)$ is an example of such a lifting and so we get a map $h^0(N) \to R \to \mathcal{O}_\mathcal{L}/\wp$ taking

- $T_q$ to $(F \circ S_p \circ \chi^{\mathrm{cyclo}}(\mathrm{Frob}_q))^{-1} \mathrm{tr}\,\rho_F(\mathrm{Frob}_q)$ for $q \nmid Np$,

- $S_q$ to $(F \circ S_p \circ \chi^{\mathrm{cyclo}}(\mathrm{Frob}_q))^{-2} \det \rho_F(\mathrm{Frob}_q)$ for $q \nmid Np$,

- $U_q$ to 0 for $q|N$,

- $U_p$ to $\phi^{-1}(\mathrm{Frob}_p)(F \circ S^p \circ \chi^{\mathrm{cyclo}})(\mathrm{Frob}_p)$, where $\phi$ is unramified and $\phi(\mathrm{Frob}_p) = F(U_p)$.

We can take $G$ to be any lifting of this map to a homomorphism $h^0(N) \to \mathcal{O}_\mathcal{L}$.

The same method of proof gives the following result.

THEOREM 3. *Let $L/\mathbb{Q}_p$ be a finite extension with ring of integers $\mathcal{O}_L$ and maximal ideal $\lambda$. Suppose that $\rho : G_\mathbb{Q} \to \mathrm{GL}_2(\mathcal{O}_L)$ is a continuous representation satisfying*

1. *$\rho$ is ramified at only finitely many primes;*

2. *$(\rho \mod \lambda)$ is irreducible;*

3. *$(\rho \mod \lambda)$ is modular;*

4. *$\rho$ is unramified at $p$ and $\rho(\mathrm{Frob}_p)$ has eigenvalues $\alpha$ and $\beta$ with distinct reductions modulo $\lambda$.*

*Then there is an integer $N$ coprime to $p$ and two homomorphisms $f_\alpha$, $f_\beta$ : $h^0(N) \to \mathcal{O}_L$ satisfying:*

1. $f_\alpha(T_q) = f_\beta(T_q) = \operatorname{tr} \rho(\operatorname{Frob}_q)$ *for all* $q \nmid Np$;

2. $f_\alpha(S_q) = f_\beta(S_q) = \det \rho(\operatorname{Frob}_q)$ *for all* $q \nmid Np$;

3. $f_\alpha(U_q) = f_\beta(U_q) = 0$ *for all* $q|N$;

4. $f_\alpha(U_p) = \alpha$ *and* $f_\beta(U_p) = \beta$.

Gross's result [Gr] (completed by [CV]) ensures the existence of $f_\alpha \bmod \lambda$ and $f_\beta \bmod \lambda$. When we argue as above and invoke Diamond's result [D] we prove the theorem.

## 2. $p$-adic modular forms

Let us first recall some facts about $p$-adic modular forms. Recall that $X_1(N)$ has a natural model over $\mathbb{Z}_p$ (in fact over $\mathbb{Z}[1/N]$). It is a natural completion of the moduli space of elliptic curves $E$ together with an embedding $i : \mu_N \hookrightarrow E[N]$. We will let $X$ denote the pull-back of $X_1(N)$ to $\mathbb{C}_p$. We will let $X_0(p)/\mathbb{C}_p$ (resp. $X(p)/\mathbb{C}_p$) denote the natural completion of the moduli space for elliptic curves $E$, an embedding $i : \mu_N \hookrightarrow E[N]$ and an isogeny $E \xrightarrow{\alpha} E'$ of degree $p$ (resp. and two points $Q_1, Q_2$ in $E[p]$ whose Weil pairing is $\zeta_p \in \mathbb{C}_p$ a fixed primitive $p^{\text{th}}$ root of 1). Note that our use of the terminology $X_0(p)$ and $X(p)$ is nonstandard, but not mentioning the level $N$ structure keeps the notation less cluttered.

There are natural maps $\pi_3 : X(p) \to X_0(p)$ (resp. $\pi_1 : X_0(p) \to X$, resp. $\pi_2 : X_0(p) \to X$) which take $(E, i, Q_1, Q_2)$ to $(E, i, E \to E/\langle Q_1 \rangle)$ (resp. $(E, i, E \xrightarrow{\alpha} E')$ to $(E, i)$, resp. $(E, i, E \xrightarrow{\alpha} E')$ to $(E', \alpha \circ i)$). These maps are all etale away from the cusps (as long as $N \geq 5$). The map $\pi_1 \circ \pi_3$ is also Galois with group $\operatorname{SL}_2(\mathbb{F}_p)$, and $\pi_3$ is thus also Galois and has group $B(\mathbb{F}_p) \subset \operatorname{SL}_2(\mathbb{F}_p)$, the subgroup of upper triangular matrices. We will let $\omega_X$ (resp. $\omega_{X_0(p)}, \omega_{X(p)}, \ldots$) denote the canonical extension to the cusps of the pull back by the identity section of the sheaf of relative differentials of the universal elliptic curve over the noncuspidal locus of $X$ (resp. $X_0(p)$, $X(p), \ldots$). Then $\pi_1^* \omega_X = \omega_{X_0(p)}$, $\pi_3^* \omega_{X_0(p)} = \omega_{X(p)}$ and there is a natural map $j = (\alpha^\vee)^* : \omega_{X_0(p)} \to \pi_2^* \omega_X$. After one inverts $p$, $j$ becomes an isomorphism. When it will not cause confusion, we shall refer to any of these sheaves as simply $\omega$.

There is a natural identification of $S_k(\Gamma_1(N)) \otimes_\mathbb{Z} \mathbb{Z}_p$ with $\Gamma(X_1(N)/\mathbb{Z}_p, \omega^{\otimes k})$. There is a map from $\mathbb{Z}_p((q))$ to $X_1(N)$ corresponding to $(\mathbb{G}_m/q^\mathbb{Z}, i^{\text{can}})/\mathbb{Z}_p((q))$. Pulling back $f \in \Gamma(X_1(N)/\mathbb{Z}_p, \omega^{\otimes k})$ to $\mathbb{Z}_p((q))$ and dividing by the generator $(dt/t)^{\otimes k}$ of $\omega^{\otimes k}/\mathbb{Z}_p((q))$ correspond to taking the $q$-expansion of $f$ at infinity. (Here $t$ is the natural parameter on $\mathbb{G}_m$ and $i^{\text{can}}$ is induced by the tautological inclusion $\mu_N \subset \mathbb{G}_m$.)

The curve $X_1(N)/\overline{\mathbb{F}}_p$ has a finite number, $s$, of supersingular points. We let $SS$ denote the union of their residue discs. This is a rigid analytic space isomorphic to the union of $s$ open discs each of radius 1. Choose parameters $T_1, ..., T_s$ such that the completed local ring of $X_1(N)/\mathbb{Z}_p^{nr}$ at the $i^{\text{th}}$ supersingular point is $\mathbb{Z}_p^{nr}[[T_i]]$. Here $\mathbb{Z}_p^{nr}$ denotes the Witt vectors of $\overline{\mathbb{F}}_p$. If $I \subseteq [0, 1)$ is a closed, open, or half-open interval with endpoints in $p^{\mathbb{Q}}$, then we define a rigid subspace $SSI \subseteq SS$ to be the union over all supersingular residue discs of the points $x$ such that if $x$ is in the $i^{\text{th}}$ disc then $|T_i(x)|_p \in I$. Here the $p$-adic absolute value is normalised so that $|p|_p = p^{-1}$. So for example, $SS[0, 1) = SS$. If $I = [0, r]$ (resp. $[0, r)$, $[r_1, r_2]$, etc.) then $SSI$ is a disjoint union of closed discs of radius $r$ (resp. open discs of radius $r$, closed annuli of radii $r_1$ and $r_2$, etc.) In general, $SSI$ will depend on the choices of the $T_i$, but if $I \subseteq (1/p, 1)$ or $[0, 1/p] \subseteq I$ then $SSI$ is independent of these choices because we have an integral model over $\mathbb{Z}_p^{nr}$. Note that we shall only consider $I$ such that $SSI$ is independent of choices in what follows. If $r \in p^{\mathbb{Q}}$ and $1 > r \geq 1/p$ we let $X_{>r} = X - SS[0, r]$. If $r \in p^{\mathbb{Q}}$ and $1 \geq r > 1/p$ we let $X_{\geq r} = X - SS[0, r)$.

We will let $E$ denote the section of $\omega_X^{\otimes(p-1)}$ over $X$ with $q$-expansion at infinity

$$1 - (2(p-1)/B_{p-1}) \sum_{n=1}^{\infty} \sigma_{p-2}(n) q^n$$

where $B_{p-1}$ denotes the Bernoulli number, and $\sigma_t(n) = \sum_{0 < d | n} d^t$. Then $E$ has a single simple zero in each connected component of $SS[0, 1/p]$ and no other zeros on $X$.

The theory of the canonical subgroup (see [Ka], particularly Theorem 3.10.7) provides rigid sections $s_1, s_2 : X_{>p^{-p/(1+p)}} \to X_0(p)$ (corresponding to $(E, i) \mapsto (E, i, E \to E/C)$ and $(E, i) \mapsto (E/C, p^{-1}i, E/C \xrightarrow{p} E)$ respectively, where $C$ denotes the canonical subgroup). These sections are isomorphisms onto their images. The induced map $\pi_1 : (s_2 X_{>p^{-1/(1+p)}}) \to X_{>p^{-p/(p+1)}}$ is finite and surjective of degree $p$. It restricts to a finite surjective degree $p$ map $s_2 X_{\geq r} \to X_{\geq r^p}$ for any $r$ with $1 \geq r > p^{-1/(p+1)}$. The induced map $\pi_1 : s_2 SS[p^{-1/(1+p)}, p^{-1/(1+p)}] \to SS[0, p^{-p/(1+p)}]$ is finite surjective of degree $p+1$. The induced map $\pi_1 : s_2 SS(p^{-p/(1+p)}, p^{-1/(1+p)}) \to SS(p^{-p/(1+p)}, p^{-1/(1+p)})$ is an isomorphism. Thus $\pi_1^{-1} SS[0, p^{-1/(1+p)}) = s_2 SS(p^{-p/(1+p)}, p^{-1/(p(1+p))})$. In fact, if $p^{-p/(1+p)} < r < p^{-1/(1+p)}$ then

$$s_2 : SS[1/(pr), r^{1/p}] \xrightarrow{\sim} \pi_1^{-1} SS[0, r].$$

We let $\mathcal{S}_k(N)$ denote the space of sections of $\omega_X^{\otimes k}$ over $X_{\geq 1}$ which vanish at each cusp. This is a $p$-adic Banach space and one of the norms in the natural equivalence class is given by $||f|| = \sup_{n > 0} |a_n(f)|$, where $f$ has $q$-expansion $\sum_{n=1}^{\infty} a_n(f) q^n$ at infinity. The operators $T(n)$ for $p \nmid n$ and the diamond operators act naturally on these spaces. One can also define an operator $U_p$ satisfying $a_n(f|U_p) = a_{np}(f)$. Hida shows that his idempotent $e$ is defined on this space and that $e\mathcal{S}_k(N)$ is finite dimensional, $U_p$ is an isomorphism on $e\mathcal{S}_k(N)$ and that $U_p$ is topologically nilpotent on $(1-e)\mathcal{S}_k(N)$. (For this and the references to Hida in the next paragraph see [H1] and [Go].)

We also let $\mathcal{S}_k^{>r}(N)$ (resp. $\mathcal{S}_k^{\geq r}(N)$) denote the space of sections of $\omega_X^{\otimes k}$ over $X_{>r}$ (resp. $X_{\geq r}$). These spaces are preserved by the Hecke operators $T(n)$ for $p \nmid n$ and by the diamond operators. For $1 \geq r > p^{-p/(1+p)}$ we define a continuous operator $V: \mathcal{S}_k^{\geq r}(N) \to \mathcal{S}_k^{\geq r^{1/p}}(N)$ as the composite

$$\Gamma(X_{\geq r}, \omega_X^{\otimes k}) \xrightarrow{\pi_1^*} \Gamma(s_2(X_{\geq r^{1/p}}), \omega_X^{\otimes k}) \xrightarrow{p^{-k}(s_2^* \circ j)} \Gamma(X_{\geq r^{1/p}}, \omega_X^{\otimes k}).$$

We have that

$$f|V = \sum_{n=1}^{\infty} a_n(f) q^{np}.$$

For $1 \geq r > p^{-p/(1+p)}$, the Hecke operator $U_p$ gives a continuous map $U_p: \mathcal{S}_k^{\geq r^{1/p}}(N) \to \mathcal{S}_k^{\geq r}(N)$. For $k > 1$ this is Corollary II.3.7 of [Go]. As Coleman explained to us, the case $k \leq 1$ can be reduced to the case $k > 1$ because the map $U_p$ is the composite:

$$\mathcal{S}_k^{\geq r^{1/p}}(N) \xrightarrow{(E^w|V)} \mathcal{S}_{k+w(p-1)}^{r^{1/p}}(N) \xrightarrow{U_p} \mathcal{S}_{k+w(p-1)}^{\geq r}(N) \xrightarrow{E^{-w}} \mathcal{S}_k^{\geq r}(N).$$

Let $S^0(N)^{(k)}$ denote the set of $F \in S^0(N)$ such that $F|S_p(x) = x^{k-1}F$ for $x \in \mu_{p-1} \subset \mathbb{Z}_p^\times$. Then Hida showed that

$$(S^0(N)^{(k)}/(u-(1+p)^{k-2})S^0(N)^{(k)}) \otimes \mathbb{C}_p \xrightarrow{\sim} e\mathcal{S}_k(N)$$

via a map taking $F$ to $\sum_{n=1}^{\infty} (F(T(n)) \bmod (u-(1+p)^{k-2}))q^n$. If $k \geq 2$ then Hida also showed that

$$e\mathcal{S}_k(N) = e(S_k(\Gamma_1(N) \cap \Gamma_0(p)) \otimes \mathbb{C}_p)$$

and hence we may deduce that

$$e\mathcal{S}_k(N) \subset \mathcal{S}_k^{>p^{-p/(p+1)}}(N).$$

The following lemma seems to be well-known to experts, but for lack of a reference we sketch a proof.

LEMMA 1. $e\mathcal{S}_k(N) \subset \mathcal{S}_k^{>p^{-p/(p+1)}}(N).$

Let $1 > r > p^{-p/(1+p)}$. The main point is to see that $U_p: \mathcal{S}_k^{\geq r}(N) \to \mathcal{S}_k^{\geq r}(N)$ is completely continuous. This follows because $\mathcal{S}_k^{\geq r}(N) \hookrightarrow \mathcal{S}_k^{\geq r^{1/p}}(N)$ is completely continuous being the composite

$$\mathcal{S}_k^{\geq r}(N) \xrightarrow{E^w} \mathcal{S}_{k+(p-1)w}^{\geq r}(N) \hookrightarrow \mathcal{S}_{k+w(p-1)}^{\geq r^{1/p}}(N) \xrightarrow{E^{-w}} \mathcal{S}_k^{\geq r^{1/p}}(N)$$

where the middle map is known to be completely continuous for $w$ sufficiently large (see [Go, Cor. I.2.9]). Moreover all eigenvalues of $U_p$ are integral, because the same is true in $\mathcal{S}_k(N)$ as follows from the $q$-expansion. Thus by Serre's theory [S1] (particularly Propositions 7 and 12), we may write $\mathcal{S}_k^{\geq r}(N) = e'\mathcal{S}_k^{\geq r}(N) \oplus (1-e')\mathcal{S}_k^{\geq r}(N)$, where $e'$ is an idempotent commuting with $U_p$, where $e'\mathcal{S}_k^{\geq r}(N)$ is finite dimensional and is spanned by generalised eigenvectors of $U_p$ with unit eigenvalues, and where $U_p$ is topologically nilpotent on $(1-e')\mathcal{S}_k^{\geq r}(N)$. Thus if $f \in \mathcal{S}_k^{\geq r}(N)$ we see that $e'f = \lim_{r \to \infty} U_p^{r!} f$ and so $e' = e|_{\mathcal{S}_k^{\geq r}(N)}$.

Let $f \in e\mathcal{S}_k(N)$. We can find a sequence $f_n \in e\mathcal{S}_{k+p^n(p-1)}(N)$ such that $f_n \to f$ (in terms of their $q$-expansions at infinity). (For example $f_n = e(fE^{p^{n+1}})$.) Then by Hida's result, for $n$ sufficiently large, $f_n \in \mathcal{S}_{k+p^n(p-1)}^{\geq r}(N)$. Finally $e(f_n/E^{p^{n+1}}) \to f$ in $e\mathcal{S}_k(N)$, but as each $e(f_n/E^{p^{n+1}})$ lies in the finite dimensional subspace $e\mathcal{S}_k^{\geq r}(N) \subset \mathcal{S}_k^{\geq r}(N)$ we also see that $f \in e\mathcal{S}_k^{\geq r}(N)$, as desired.

We will now state and prove our key technical result.

THEOREM 4. *Suppose $k \in \mathbb{Z}$, $N \in \mathbb{Z}_{\geq 5}$ and $p \geq 5$ is a prime not dividing $N$. Suppose $\alpha$ and $\beta$ are distinct nonzero elements of $\mathbb{C}_p$. Suppose that $f_\alpha, f_\beta \in \mathcal{S}_k^{>t}(N)$ for some $t < 1$ are eigenvectors for $U_p$ with eigenvalues $\alpha$ and $\beta$. Suppose finally that for all positive integers $n$ not divisible by $p$,*

$$a_n(f_\alpha) = a_n(f_\beta).$$

*Then $f = (\alpha f_\alpha - \beta f_\beta)/(\alpha - \beta)$ is classical, i.e. lies in $S_k(\Gamma_1(N)) \otimes \mathbb{C}_p$.*

We first note that we may take $t = p^{-p/(1+p)}$ (because $\alpha$ and $\beta$ are nonzero). Set $f' = (f_\alpha - f_\beta)/(\alpha - \beta)$ so that $f' = f|V \in \mathcal{S}_k^{>p^{-1/(1+p)}}(N)$.

Choose $r, r' \in p^\mathbb{Q}$ with

$$p^{-p/(1+p)} < r' < r < p^{-1/(1+p)}.$$

Let $S$ denote $(\pi_1 \circ \pi_3)^{-1} SS[0,r] \subset X(p)$. Define a section $g$ of $\omega_{X(p)}^{\otimes k}$ over $S$ by $g = \pi_3^* \circ j^{-1} \circ \pi_2^*(p^k f')$. We will show below that $g$ is invariant for the action of $\mathrm{SL}_2(\mathbb{F}_p)$. Because $S \to SS[0,r]$ is finite it will then follow that $g = (\pi_3 \circ \pi_1)^*(h)$ for some section $h$ of $\omega_X$ on $SS[0,r]$. On the other hand, on $s_2(SS(p^{-1/(1+p)}, r^{1/p}])$ we have $\pi_2^*(p^k f') = j \circ \pi_1^*(f)$ and hence on $\pi_3^{-1} s_2(SS(p^{-1/(1+p)}, r^{1/p}])$ we have $g = (\pi_3 \circ \pi_1)^*(f)$. Thus $h|_{SS(p^{-p/(1+p)}, r]} = f$. Then we glue the sections $h$ on $SS[0,r]$ and $f$ on $X_{\geq r'}$ to give the desired section $f \in S_k(\Gamma_1(N)) \otimes_\mathbb{Z} \mathbb{C}_p$ (by rigid GAGA, [Ko]).

We now turn to the proof that $g$ is invariant under $\mathrm{SL}_2(\mathbb{F}_p)$. If $C$ is a connected component of $S$, let $G_C$ denote the subgroup of $\mathrm{SL}_2(\mathbb{F}_p)$ stabilising $C$. Note that $B(\mathbb{F}_p)$ acts transitively on the connected components of $S$ lying over any given connected component of $SS[0,r]$ and so $\mathrm{SL}_2(\mathbb{F}_p) = B(\mathbb{F}_p).G_C$.

On the other hand $B(\mathbb{F}_p)$ fixes $g$. Thus it will suffice to show that $G_C$ fixes $g|_C$ for any connected component $C$. Let $\sigma \in G_C$ and let

$$C' = \pi_3^{-1} s_2(SS(p^{-1/(1+p)}, r^{1/p}]) \cap \sigma^{-1} \pi_3^{-1} s_2(SS(p^{-1/(1+p)}, r^{1/p}]).$$

On $C'$ we have

$$\sigma^*(g) = \sigma^* \pi_3^* j^{-1} \pi_2^*(p^k f') = \sigma^* \pi_3^* \pi_1^*(f) = \pi_3^* \pi_1^*(f) = \pi_3^* j^{-1} \pi_2^*(p^k f') = g.$$

Thus it suffices to show that $C'$ is infinite. Suppose that $\pi_3 : C \to \pi_1^{-1} SS[0, r]$ has degree $d$. Then

$$\pi_1 \pi_3 : \pi_3^{-1} s_2(SS(p^{-1/(1+p)}, r^{1/p}]) \to SS(p^{-p/(1+p)}, r]$$

is of degree $dp$, as is

$$\pi_1 \pi_3 : \sigma^{-1} \pi_3^{-1} s_2(SS(p^{-1/(1+p)}, r^{1/p}]) \to SS(p^{-p/(1+p)}, r].$$

Hence, since the degree of $C \to SS[0, r]$ is $d(p+1)$, it follows by an elementary argument that the fibres of $\pi_1 \pi_3 : C' \to SS(p^{-p/(1+p)}, r]$ have cardinality at least $d(p-1)$ and $C'$ surjects onto $SS(p^{-p/(1+p)}, r]$. The theorem follows.

Before proceeding to the proof of our main theorem we note that Theorem 4 gives a partial answer to a question of Gouvea. We will call two normalised eigenforms $f_1, f_2 \in \mathcal{S}_k(N)$ equivalent away from $p$ if $a_n(f_1) = a_n(f_2)$ for all $n$ not divisible by $p$. Gouvea notes ([Go, §II.3.3]) that if $f$ is any normalised eigenform in $\mathcal{S}_k(N)$ and if $\alpha \in \mathbb{C}_p$ satisfies $|\alpha|_p < 1$ then there is an eigenform $f_\alpha \in \mathcal{S}_k(N)$ equivalent away from $p$ to $f$ with $f_\alpha | U_p = \alpha f_\alpha$. Gouvea asks how many of these eigenforms $f_\alpha$ can be overconvergent, i.e. lie in $\mathcal{S}_k^{>r}(N)$ for some $r < 1$. Our theorem implies the following result.

COROLLARY 1. *With the notation as above, and for fixed $f$, $f_\alpha$ can be overconvergent for at most two nonzero values of $\alpha$. If $f_\alpha$ and $f_\beta$ are overconvergent for two distinct and nonzero $\alpha$ and $\beta$ then $(\alpha f_\alpha - \beta f_\beta)/(\alpha - \beta) \in S_k(\Gamma_1(N)) \otimes \mathbb{C}_p$. In particular $k \geq 1$, $\alpha\beta = p^{k-1}\chi(p)$ where $f|\langle p \rangle = \chi(p)f$, and the Galois representation $\rho_f$ associated to $f$ is crystalline at $p$.*

## 3. Weight one forms

Putting the results of the last two sections together we are now in a position to prove our main theorem. Let $L/\mathbb{Q}_p$ be a finite extension with ring of integers $\mathcal{O}_L$ and maximal ideal $\lambda$.

THEOREM 5. *Suppose that $\rho : G_\mathbb{Q} \to \mathrm{GL}_2(\mathcal{O}_L)$ is a continuous representation satisfying the following conditions.*

1. *$\rho$ ramifies at only finitely many primes.*

2. $(\rho \bmod \lambda)$ *is modular and irreducible.*

3. $\rho$ *is unramified at $p$ and $\rho(\mathrm{Frob}_p)$ has eigenvalues $\alpha$ and $\beta$ with distinct reductions modulo $\lambda$.*

*Then there exists an integer $N$ coprime to $p$ and an eigenform $f \in S_1(\Gamma_1(N)) \otimes_{\mathbb{Z}} \mathcal{O}_L$ such that for almost all primes $q$, $a_q(f) = \operatorname{tr} \rho(\mathrm{Frob}_q)$. In particular $\rho$ has finite image and for any embedding $i : L \hookrightarrow \mathbb{C}$ the Artin L-function $L(i \circ \rho, s)$ is entire.*

Combining Theorem 3 and Lemma 1 we can find such an integer $N$ and two sections $f_\alpha$ and $f_\beta$ in $e\mathcal{S}_1^{>-p/(p+1)}(N)$ which are eigenvectors for the Hecke operators $T(n)$ for $p \nmid n$ and for $U_p$ and which have the following eigenvalues.

- $f_\alpha | T_q = (\operatorname{tr} \rho(\mathrm{Frob}_q)) f_\alpha$ and $f_\beta | T_q = (\operatorname{tr} \rho(\mathrm{Frob}_q)) f_\beta$ if $q \nmid Np$.

- $f_\alpha | \langle q \rangle = (\det \rho(\mathrm{Frob}_q)) f_\alpha$ and $f_\beta | \langle q \rangle = (\det \rho(\mathrm{Frob}_q)) f_\beta$ if $q \nmid Np$.

- $f_\alpha | U_q = 0$ and $f_\beta | U_q = 0$ if $q | N$.

- $f_\alpha | U_p = \alpha f_\alpha$ and $f_\beta | U_p = \beta f_\beta$.

Then by Theorem 4 we see that

$$f = (\alpha f_\alpha - \beta f_\beta)/(\alpha - \beta) \in S_1(\Gamma_1(N)) \otimes_{\mathbb{Z}} \mathcal{O}_L$$

is the desired form.

COROLLARY 2. *Let $\mathcal{L}$ be a finite extension of the fraction field of $\Lambda$ and let $\wp$ be a height one prime of $\mathcal{O}_\mathcal{L}$ above $((1+p)u - 1)$. Suppose that $F \in S^0(N) \otimes \mathcal{L}$ is a non-Eisenstein $\Lambda$-adic eigenform for which $F|\langle x \rangle = x^{-1} F$ for $x \in \mu_{p-1} \subset \mathbb{Z}_p^\times$. Suppose also that $(\rho_F \bmod \wp)$ is unramified at $p$ and the eigenvalues of $(\rho_F \bmod \wp)(\mathrm{Frob}_p)$ are distinct modulo the maximal ideal of $\mathcal{O}_\mathcal{L}/\wp$. Then $\sum_{n=1}^\infty (F(T(n)) \bmod \wp) q^n \in S_1(\Gamma_1(Np)) \otimes_{\mathbb{Z}} (\mathcal{O}_\mathcal{L}/\wp)$.*

Gross has pointed out to us the following consequence of our main theorem. It partially answers a question posed to Gross by Serre.

COROLLARY 3. *Suppose that $p > 5$ is a prime and $N \geq 5$ is an integer not divisible by $p$. Suppose that $\overline{f} \in H^0(X_1(N) \times_{\mathbb{Z}_p} \overline{\mathbb{F}}_p, \omega)$ is a normalised eigenform such that $\overline{f}|U_r = 0$ for all primes $r | N$ and such that $X^2 - a_p(\overline{f}) X + \chi(p)$ has distinct roots, where $\overline{f}|\langle p \rangle = \chi(p) \overline{f}$. Suppose also that the Galois representation $\overline{\rho}_{\overline{f}} : G_{\mathbb{Q}} \to \mathrm{GL}_2(\overline{\mathbb{F}}_p)$ associated to $\overline{f}$ is irreducible. Then $\overline{f}$ is in the image of*

$$S_1(\Gamma_1(N)) \otimes_{\mathbb{Z}} \overline{\mathbb{F}}_p \hookrightarrow H^0(X_1(N) \times_{\mathbb{Z}_p} \overline{\mathbb{F}}_p, \omega)$$

*if and only if $p$ does not divide the order of $\overline{\rho}_{\overline{f}}(G_{\mathbb{Q}})$.*

(Before explaining the proof we remark that the condition that $\overline{f}|U_r = 0$ for $r|N$ can almost certainly be suppressed, but we have not been through the slightly tedious details. We also remark that we could prove a similar theorem when $p = 5$, but in that case the conclusion would have to be slightly modified. We leave the details to the reader.)

First suppose that $\overline{f}$ is in the image of

$$S_1(\Gamma_1(N)) \otimes_{\mathbb{Z}} \overline{\mathbb{F}}_p \hookrightarrow H^0(X_1(N) \times_{\mathbb{Z}_p} \overline{\mathbb{F}}_p, \omega).$$

Then we can find a finite extension $L/\mathbb{Q}_p$ and a normalised eigenform $f \in S_1(\Gamma_1(N)) \otimes \mathcal{O}_L$ which reduces to $\overline{f}$ (by the "going up theorem" applied to the Hecke algebra). The Galois representation $\rho_f : G_{\mathbb{Q}} \to \mathrm{GL}_2(\mathcal{O}_L)$ associated to $f$ by Deligne and Serre [DS] has finite image and reduces to $\overline{\rho}_{\overline{f}}$. By the classification of finite subgroups of $\mathrm{GL}_2(\mathbb{C})$ we see that the image of $\rho_f$ can not have $\mathrm{PSL}_2(\mathbb{F}_p)$ as a subquotient (as $p > 5$), and hence the image of $\overline{\rho}_{\overline{f}}$ can not contain any conjugate of $\mathrm{SL}_2(\mathbb{F}_p)$. The classification of finite subgroups of $\mathrm{GL}_2(\overline{\mathbb{F}}_p)$ then shows that either $\overline{\rho}_{\overline{f}}$ is reducible or $p$ does not divide the order of $\overline{\rho}_{\overline{f}}(G_{\mathbb{Q}})$.

Conversely suppose that $p$ does not divide the order of $\overline{\rho}_{\overline{f}}(G_{\mathbb{Q}})$. Then a simple cohomological calculation shows that we can find a finite extension $L/\mathbb{Q}_p$ and a lifting $\rho : G_{\mathbb{Q}} \to \mathrm{GL}_2(\mathcal{O}_L)$ of $\overline{\rho}_{\overline{f}}$ such that the image of $\rho$ maps isomorphically to the image of $\overline{\rho}_{\overline{f}}$. (If $\overline{\rho}_{\overline{f}}$ is defined over $\mathbb{F}_q$ then any $L$ containing the Witt vectors of $\mathbb{F}_q$ will do.) By Proposition 2.7 of [E], $\overline{\rho}_{\overline{f}}$ is unramified at $p$ and $\overline{\rho}(\mathrm{Frob}_p)$ has two distinct eigenvalues, namely the two roots of $X^2 - a_p(\overline{f})X + \chi(p)$. Thus $\rho$ is also unramified at $p$, and so by our main theorem we can find a normalised eigenform $f \in S_1(\Gamma_1(M)) \otimes \mathcal{O}_L$ for some positive integer $M$, such that $a_r(f) = \mathrm{tr}\,\rho(\mathrm{Frob}_f)$ for all but finitely many primes $r$. By the theory of newforms we may in fact take $M$ equal to the conductor of $\rho$ and ensure that $a_r(f) = \mathrm{tr}\,\rho^{I_r}(\mathrm{Frob}_r)$ for all primes $r$. Alternatively we may take $M$ to be the product of the conductor of $\rho$ and the product over all primes $r|N$ of $r^{\dim \rho^{I_r}}$. We thus ensure that $a_r(f) = \mathrm{tr}\,\rho^{I_r}(\mathrm{Frob}_r)$ if $r \nmid N$, while $a_r(f) = 0$ if $r|N$. On the other hand because $\overline{f}|U_r = 0$ if $r|N$ we have that the product of the conductor of $\overline{\rho}_{\overline{f}}$ and $\prod_{r|N} r^{\dim \overline{\rho}_{\overline{f}}^{I_r}}$ divides $N$. Also $\rho$ and $\overline{\rho}_{\overline{f}}$ have the same conductor (as $p$ does not divide the order of $\rho(G_{\mathbb{Q}})$ and $\rho$ is unramified at $p$). Thus this last form $f$ has level dividing $N$ and reduces modulo $p$ to $\overline{f}$.

University of California at Berkeley, Berkeley, CA
*Current address*: Imperial College, London, England
*E-mail address*: buzzard@ic.ac.uk

Harvard University, Cambridge, MA
*E-mail address*: rtaylor@math.harvard.edu